\documentclass{article}
\usepackage{amsmath}
\usepackage{amsfonts}
\usepackage{amssymb}
\usepackage{url}
\usepackage{color}
\usepackage{enumerate}
\usepackage{graphics}
\usepackage{graphicx}
\usepackage{comment}

\newtheorem{theorem}{Theorem}[section]
\newtheorem{lemma}[theorem]{Lemma}

\newtheorem{definition}[theorem]{Definition}
\newtheorem{corollary}[theorem]{Corollary}

\newcommand{\NP}{{\mathcal N}{\mathcal P}}

\newcommand{\calE}{{\mathcal E}}
\newcommand{\Z}{{\mathbb Z}}

\newcommand{\X}{{\mathbb X}}

\newcommand{\R}{{\mathbb R}}
\def\1{{\mathbf 1}}

\DeclareMathOperator{\pa}{pa}

\DeclareMathOperator{\portrait}{portrait}
\DeclareMathOperator{\pat}{pat}

\numberwithin{equation}{section}

\def\Orthant_j{{\mathcal O}_{j}}

\def\DAG{{\mbox{DAG}}}

\def\chv{{\ve c}}
\def\staim{u}

\def\vex{{\ve x}}

\usepackage{enumerate}
\def\ve#1{\mathchoice{\mbox{\boldmath$\displaystyle\bf#1$}}
{\mbox{\boldmath$\textstyle\bf#1$}}
{\mbox{\boldmath$\scriptstyle\bf#1$}}
{\mbox{\boldmath$\scriptscriptstyle\bf#1$}}}

\newcommand{\boproof}{\textbf{Proof.} }
\newcommand{\eoproof}{\hspace*{\fill} $\square$ \vspace{5pt}}

\begin{document}
\setlength{\parindent}{0pt} \setlength{\parskip}{2ex plus 0.4ex
minus 0.4ex}

\title{Learning restricted Bayesian network structures}
\author{Raymond~Hemmecke (TU Munich, hemmecke@ma.tum.de)\\ Silvia~Lindner (TU Munich, slindner@ma.tum.de)\\ Milan~Studen\'{y} (UTIA Prague, studeny@utia.cas.cz)}

\date{\today}

\maketitle

\begin{abstract}
  Bayesian networks are basic graphical models, used widely both in statistics and artificial intelligence. These statistical models of conditional independence structure are described by acyclic directed graphs whose nodes correspond to (random) variables in consideration. A quite important topic is the learning of Bayesian network structures, which is determining the best fitting statistical model on the basis of given data. Although there are learning methods based on statistical conditional independence tests, contemporary methods are mainly based on maximization of a suitable quality criterion that evaluates how good the graph explains the occurrence of the observed data. This leads to a nonlinear combinatorial optimization problem that is in general $\NP$-hard to solve.

  In this paper we deal with the complexity of learning restricted Bayesian network structures, that is, we wish to find network structures of highest score within a given subset of all possible network structures. For this, we introduce a new unique algebraic representative for these structures, called the characteristic imset. We show that these imsets are always $0$-$1$-vectors and that they have many nice properties that allow us to simplify long proofs for some known results and to easily establish new complexity results for learning restricted Bayes network structures.
\end{abstract}

\section{Introduction}

  {\em Bayesian networks\/} are basic graphical models, used widely both in statistics \cite{bib:Lau96} and artificial intelligence \cite{bib:Pea88}. These statistical models of conditional independence structure are described by acyclic directed graphs whose nodes correspond to (random) variables in consideration.

  A quite important topic is {\em learning Bayesian network structures\/} \cite{bib:Nea04}, which is determining the statistical model on the basis of given data. Although there are learning methods based on statistical conditional independence tests, contemporary methods are mainly based on maximization of a suitable {\em quality criterion} or {\em score function} $\mathcal{Q}(G,D)$ of the (acyclic directed) graph $G$ and the (given = fixed) data $D$, evaluating how good the graph $G$ explains the occurrence of the observed data $D$. This leads to a nonlinear combinatorial optimization problem that is $\NP$-hard \cite{bib:Chi96,bib:CHM04}. Below we will consider learning restricted Bayesian network structures. Some of these problems remain $\NP$-hard while others are polynomial-time solvable.

  It may happen that two different acyclic directed graphs describe the same statistical model, that is, they are {\em Markov equivalent}. A classic result \cite{bib:Fry90,bib:VerPea91} says that two acyclic directed graphs are Markov equivalent if and only if they have the same underlying undirected graph and the same set of immoralities (= special induced subgraphs $a\to c\leftarrow b$ over three nodes $\{a,b,c\}$ with no arc between $a$ and $b$ in either direction). In order to remove this ambiguity of Markov equivalent models/graphs, one is interested in having a unique representative for each Bayesian network structure (= statistical model). A classic unique graphical representative is the {\em essential graph\/} \cite{bib:AMP97} of the corresponding Markov equivalence class of acyclic directed graphs, which is a special graph allowing both directed and undirected edges (see Section \ref{Ssec.Learning} for more details).

  Any reasonable score function should be {\em score-equivalent} \cite{bib:Bou95}, that is, $\mathcal{Q}(G,D)=\mathcal{Q}(H,D)$ for any two Markov equivalent graphs $G$ and $H$. Another standard technical requirement is that the criterion has to be (additively) {\em decomposable} into contributions from the parent sets $\pa_G(i)$ of each node $i$ \cite{bib:Chi02} (see Section \ref{Ssec.Learning} for more details).

  In this paper, we deal with learning (restricted) decomposable models \cite{bib:Lau96}, interpreted as Bayesian network structures. Decomposable models are exactly those models whose essential graph is an \emph{undirected} (and thus also necessarily {\em chordal}) graph. That is, decomposable models correspond to graphical models without immoralities. As input to our learning problem we assume that we are given an undirected graph $K$ and an evaluation oracle for the score function $\mathcal{Q}(\cdot,D)$. Note that we do not assume the actual data $D$ being part of the input itself. Of course, the evaluation oracle uses the given data $D$ in order to evaluate score function values. However, in our treatment, we remove the complexity of evaluating score function values from the overall complexity. In particular, this means that the (large or small) number of data vectors in $D$ will be irrelevant for our complexity results.

  We show that learning spanning trees of $K$ and learning forests in $K$ are both polynomial-time solvable. For learning spanning trees of $K$, this observation has been already made in \cite{bib:ChowLiu68} for {\em specific} score functions. Moreover, we show that if we impose degree bounds $\deg(v)\leq k$ on all nodes $v\in N$ for some constant $k\geq 2$, then both problems become $\NP$-hard. We also show that learning chordal subgraphs of $K$ is $\NP$-hard. This result, however, has been already shown even for specific score functions and also for the case of fixed bounded size of possible cliques \cite{bib:Sre01}. We include our short proof to emphasize the simplicity and usefulness of our approach to easily recover also hardness results.

  We will rewrite the nonlinear combinatorial optimization problem behind the learning problem into a linear integer optimization problem (in higher dimension) by using an algebraic approach to the description of conditional independence structures \cite{bib:Stud05} that represents them by certain vectors with integer components, called {\em imsets} (short for ``Integer Multi-SETs''). In the context of learning Bayesian networks this led to the proposal to represent each Bayesian network structure  uniquely by a so-called {\em standard imset}. The advantage of this algebraic approach is that every reasonable score function (score equivalent and decomposable), becomes an affine function of the standard imset (see Chapter 8 in \cite{bib:Stud05}). Moreover, it has recently been shown in \cite{bib:StudVomHem10} that the standard imsets over a fixed set of variables are exactly the vertices of their convex hull, the {\em standard imset polytope}. These results allow one to apply the methods of polyhedral geometry in the area of learning Bayesian networks, because they transform this task to a linear programming problem.

  Instead of considering standard imsets, we introduce a different unique representative that is obtained from the standard imset by an invertible affine linear map that preserves lattice points in both directions. We call these new representatives \emph{characteristic imsets}, as they are $0$-$1$-vectors and as they also contain, for each acyclic directed graph, the characteristic vector of the underlying undirected graph. Although, mathematically, this map is simply a change in coordinates, the characteristic imset is much closer to the graphical description because it allows one to identify immediately both the underlying undirected graph and the immoralities. Our procedure for recovering the essential graph from the characteristic imset is much simpler than the reconstruction from the standard imset as presented in \cite{bib:StudVom09}.

  Moreover, due to the affine transformation, every reasonable score function is also an affine function of the characteristic imset. Thus, learning Bayesian network structures can be reduced to solving a linear optimization problem over a certain $0$-$1$-polytope. Unfortunately, a complete facet description for this polytope (for general $|N|$) is still unknown. A conjectured list of all facets for the standard imset polytope (and consequently also for the characteristic imset polytope) is presented in \cite{bib:StuVom11}. A complete facet description is also unknown for the convex hull of all characteristic imsets of undirected chordal graphs, although the characteristic imsets themselves are well-understood in this case (see Section \ref{sec.characteristic}).

  To summarize, we offer a new method for analyzing the learning procedure through an algebraic way of representing statistical models. We believe that our approach via characteristic imsets brings a tremendous mathematical simplification that allows us to easily recover known results and to establish new complexity results. We also think that a better understanding of the polyhedral properties of the characteristic imset polytope (complete facet description or all edge directions) will lead to future applications of efficient (integer) linear programming methods and software in this area of learning Bayesian network structures.

\section{Basic concepts}\label{sec.basic}

  We tacitly assume that the reader is familiar with basic concepts from polyhedral geometry. We only recall briefly the definitions of concepts mentioned above, but skip their statistical motivation.

  Throughout the paper $N$ is a finite non-empty set of {\em variables}; to avoid the trivial case we assume $|N|\geq 2$. In statistical context, the elements of $N$ correspond to random variables in consideration; in graphical context, they correspond to nodes.

\subsection{Graphical concepts}

  Graphs considered here have a finite non-empty set of nodes $N$ and two types of edges: directed edges, called {\em arcs} (or {\em arrows} in machine learning literature), denoted by $i\rightarrow j$ (or $j\leftarrow i$), and {\em undirected edges}. No loops or multiple edges between two nodes are allowed.

  A set of nodes $C\subseteq N$ is a {\em clique} (or a {\em complete set}) in $G$ if every pair of distinct nodes in $C$ is connected by an undirected edge. An {\em immorality} in a graph $G$ is an induced subgraph (of $G$) for three nodes $\{ a,b,c\}$ in which $a\rightarrow c\leftarrow b$ and $a$ and $b$ are not adjacent. An undirected graph is called {\em chordal}, if every (undirected) cycle of length at least $4$ has a chord, that is, an edge connecting two non-consecutive nodes in the cycle. A {\em forest} is an undirected graph without undirected cycles. A connected forest over $N$ is called a {\em spanning tree}. By the {\em degree $\deg_G(i)$} of a node $i\in N$ (in an undirected graph $G$), we mean the number of edges incident to $i$ in $G$.

  Note that an undirected graph is chordal if and only if all its edges can be directed in such a way that the result is an acyclic directed graph without immoralities (see \S\,2.1 in \cite{bib:Lau96}).

  Occasionally, we will use the (in the machine learning community) commonly used acronym ``\DAG'' for ``directed acyclic graph'', although the grammatically correct phrase is ``acyclic directed graph".

\subsection{Learning Bayesian network structures}\label{Ssec.Learning}

  In statistical context, to each variable (= node) $i\in N$ is assigned a finite (individual) sample space $\X_{i}$ (= the set of possible values); to avoid technical problems assume $|\X_{i}|\geq 2$, for each $i\in N$. A {\em Bayesian network structure\/} defined by a \DAG\ $G$ (over $N$) is formally the class of discrete probability distributions $P$ on the joint sample space $\prod_{i\in N} \X_{i}$ that are Markovian with respect to $G$. Note that $P$ is {\em Markovian\/} with respect to $G$ if it satisfies conditional independence restrictions determined by the respective separation criterion (see \cite{bib:Lau96, bib:Pea88}).

  Different \DAG s over $N$ can be {\em Markov equivalent}, which means they define the same Bayesian network structure. The classic graphical characterization of (Markov) equivalent graphs is this: they are equivalent if and only if they have the same underlying undirected graph and the same immoralities (see \cite{bib:AMP97}). The classic unique graphical representative of a Bayesian network structure is the {\em essential graph} $G^{*}$ of the respective (Markov) equivalence class $\mathcal{G}$ of acyclic directed graphs: one has $a\to b$ in $G^{*}$ if this arc occurs in every graph from $\mathcal{G}$ and it has an undirected edge between $a$ and $b$ in $G^{*}$ if one has $a\rightarrow b$ in one graph and $b\rightarrow a$ in another graph (from $\mathcal{G}$).
  A less informative (unique) representative is the {\em pattern} $\pat(G)$ (of any $G$ in $\mathcal{G}$): it is obtained from the underlying graph of $G$ by directing (only) those edges that belong to immoralities (in $G$).

  {\em Learning a Bayesian network structure\/} means to determine it on the basis of an observed (complete) database $D$ (of length $\ell\geq 1$), which is a sequence $x_{1},\ldots ,x_{\ell}$ of elements of the joint sample space. $D$ is called {\em complete} if all components of the elements $x_{1},\ldots,x_{\ell}$ are known. A {\em quality criterion} is a real function $\mathcal{Q}$ of two variables: of an acyclic directed graph $G$ and of a database $D$. A learning procedure consists in maximizing the function $G\mapsto \mathcal{Q}(G,D)$ for given fixed $D$. Since the aim is to learn a Bayesian network structure, the criterion should be {\em score equivalent}, which means, $\mathcal{Q}(G,D)=\mathcal{Q}(H,D)$ for any pair of Markov equivalent graphs $G,H$ and for any database $D$. A standard technical requirement \cite{bib:Chi02} is that the criterion has to be (additively) {\em decomposable}, which means, it can be written as follows:
  \[
    \mathcal{Q} (G,D) =\sum_{i\in N} q_{i|\pa_{G}(i)} (D_{\{ i\}\cup\pa_{G}(i)}),
  \]
  where $D_{A}$ denotes the projection of the database $D$ to $\prod_{i\in A} \X_{i}$ (for $\emptyset\neq A\subseteq N$) and $q_{i|B}$ for $i\in N$, $B\subseteq N\setminus\{ i\}$ are real functions.

  Finally, let us remark that the essential graph $G^*$ of a \DAG\, $G$ is an undirected graph if and only if $G$ has no immoralities. Consequently, every cycle in the undirected graph underlying $G$ (which must coincide with $G^*$) of length at least $4$ must contain a chord (otherwise there exists an immorality on this cycle in $G$). Therefore, if an essential graph is undirected it has to be \emph{chordal}. Conversely, if $G^*$ is chordal, $G$ cannot have an immorality. Therefore, {\em learning decomposable models\/} can be viewed as learning (special) Bayesian network structures corresponding to chordal undirected essential graphs \cite{bib:AMP97b}.

\subsection{Algebraic approach to learning}\label{Subsection: Algebraic approach to learning}

  An {\em imset} over $N$ is a vector in $\Z^{2^{|N|}}$, whose components are indexed by subsets of $N$. Traditionally, all subsets of $N$ are considered, although in Section \ref{sec.characteristic} we also consider imsets with a restricted domain (components corresponding to the empty set and to singletons are dropped, since they linearly depend on the other components). Every vector in $\R^{2^{|N|}}$ can be written as a (real) combination of basic vectors $\delta_{A}\in \{0,1\}^{2^{|N|}}$:
  \[
    \delta_{A}(T) = \left\{
    \begin{array}{ll}
      1 & ~~\mbox{\rm if}~~ T=A\,,\\
      0 & ~~\mbox{\rm if}~~ T\subseteq N,\; T\neq A\,,
    \end{array} \right.
    \quad \mbox{for $T\subseteq N$ (if $A\subseteq N$ is fixed).}
  \]
  This allows us to give formulas for imsets. Given an acyclic directed graph $G$ over $N$, the {\em standard imset} for $G$ is given by
  \begin{equation}
    \staim_{G}=\delta_{N}-\delta_{\emptyset} + \sum_{i\in N}
    \left\{\,
      \delta_{\pa_{G}(i)} - \delta_{\{ i\}\cup\pa_{G}(i)}\,
    \right\},
    \label{eq.staim}
  \end{equation}
  where the basic vectors can cancel each other. It is a unique algebraic representative of the corresponding Bayesian network structure because $\staim_{G}=\staim_{H}$ if and only if $G$ and $H$ are Markov equivalent (Corollary 7.1 in \cite{bib:Stud05}). The convex hull of the set of all standard imsets over $N$ is the {\em standard imset polytope}.

  An important result from the point of view of an algebraic approach to learning Bayesian network structures is that any score equivalent and decomposable quality criterion (= score function) $\mathcal{Q}$ has the form
  \begin{equation}
    \mathcal{Q}(G,D) \,=\,
      s^{\mathcal{Q}}_{D} -\langle t^{\mathcal{Q}}_{D},\staim_{G}\rangle\,,
    \label{eq.criter}
  \end{equation}
  where $\langle \ast,\ast\rangle$ denotes the scalar product, and both $s^{\mathcal{Q}}_{D}\in\R$ and $t^{\mathcal{Q}}_{D}\in \R^{2^{|N|}}$ only depend on the database $D$ and the chosen quality criterion (see Lemmas 8.3 and 8.7 in \cite{bib:Stud05}). In particular, the task to maximize $\mathcal{Q}$ is equivalent to finding the optimum of a linear function over the standard imset polytope.

\section{Characteristic imsets}\label{sec.characteristic}

  In this section we introduce the notion of a \emph{characteristic imset\/} and prove some useful facts about it. For example, we show that this imset is always a $0$-$1$ vector.

  \begin{definition}\rm\label{def.char-im}
    Given an acyclic directed graph $G$ over $N$, let $\staim_{G}$ be the standard imset for $G$. We introduce
    \begin{eqnarray*}
      \portrait[\staim_{G}] &:=& (\,\portrait[\staim_{G}]\,(T)\,)_{T\subseteq N, |T|>1}\;\;\;\in\Z^{2^{|N|}-|N|-1}, ~\mbox{with}\\
      \portrait[\staim_{G}]\,(T) &:=& \sum_{X\subseteq N: T\subseteq X} \staim_{G}(X) ~~\mbox{for~} T\subseteq N,\ |T|>1,
    \end{eqnarray*}
    and call $\portrait[\staim_{G}]$ the \emph{upper portrait} of $\staim_{G}$ or, simply, of $G$.

    Moreover, we will call
    \[
      \chv_{G}:=\ve 1-\portrait[\staim_{G}]\in\Z^{2^{|N|}-|N|-1}
    \]
    the \emph{characteristic imset} of $G$.
  \end{definition}

  Characteristic imsets are unique representatives of Markov equivalence classes. This is because the standard imset are unique representatives and because the upper portrait map is an affine linear map that is invertible. The inverse map is given by the well-known Mobius inversion formula \cite{Bender+Goldman}. In fact, both maps assign lattice points to lattice points!

  Characteristic imsets have remarkable properties and, as we will show below, their entries directly encode the underlying undirected graph and the immoralities of the given acyclic directed graph.

  \begin{theorem}\label{Theorem: Characteristic imsets are 0-1}
    Let $G$ be an acyclic directed graph over $N$. For any $T\subseteq N$, $|T|>1$ we have $\chv_{G}(T)\in \{0,1\}$ and $\chv_{G}(T)=1$ iff there exists some $i\in T$ with $T\setminus \{ i\}\subseteq \pa_G(i)$. In particular, $\chv_{G}\in \{0,1\}^{2^{|N|}-|N|-1}$.
  \end{theorem}

  \boproof Consider the defining formula (\ref{eq.staim}) for the standard imset. For any $T\subseteq N$, $|T|>1$, the value $\portrait[\staim_{G}]\,(T)$ can be computed as
  \[
    \portrait[\staim_{G}]\,(T)=\sum_{X\subseteq N: T\subseteq X} \staim_{G}(X)=1+\sum_{i\in N:T\subseteq\pa(i)} 1 - \sum_{i\in N:T\subseteq\pa(i)\cup\{i\}} 1.
  \]
  Hence, we get
  \begin{eqnarray*}
    \chv_{G}(T) & = & 1- \portrait[\staim_{G}]\,(T)\\
    & = & \sum_{i\in N:T\subseteq\pa(i)\cup\{i\}} 1 - \sum_{i\in N:T\subseteq\pa(i)} 1\\
    & = & \sum_{i\in N:T\subseteq\pa(i)\cup\{i\},i\in T} 1\\\
    & = & \sum_{i\in T:T\setminus\{ i\}\subseteq\pa(i)} 1.
  \end{eqnarray*}
  For fixed $T$, assume that there are two different elements $i,j\in T$ with $T\setminus\{ i\}\subseteq \pa_G(i)$ and $T\setminus\{j\}\subseteq\pa_G(j)$. This implies both $i\in\pa_G(j)$ and $j\in\pa_G(i)$. The simultaneous existence of the arcs $i\rightarrow j$ and $j\to i$, however, contradicts the assumption that $G$ is acyclic. Thus, for each $T\subseteq N$, there is at most one $i\in T$ with $T\setminus\{ i\}\subseteq \pa_G(i)$. Consequently,
  \[
    \chv_{G}(T)=\sum_{i\in T:T\setminus\{ i\}\subseteq\pa(i)} 1\;\;\; \in\{0,1\},
  \]
  and thus $\chv_{G}\in\{0,1\}^{2^{|N|}-|N|-1}$. \eoproof

  \begin{corollary}\label{Cor: Vertices are the only lattice points}
    For any $N$, the only lattice points in the standard imset polytope and in the characteristic imset polytope are their vertices.
  \end{corollary}

  \boproof The statement holds for any $0$-$1$-polytope and thus in particular also for the characteristic imset polytope. Moreover, the portrait map and its inverse, the Mobius map, are affine linear maps between $\staim_{G}$ and $\chv_{G}$ that map lattice points to lattice points. Thus, the result holds also for the standard imset polytope. \eoproof

  {\bf Remark.} Corollary \ref{Cor: Vertices are the only lattice points} (for standard imsets) has already been stated and proved in \cite{bib:StuVom11}. The original proof of this result in the manuscript of \cite{bib:StuVom11} was quite long and complicated. Later discussions among the authors of the present paper led to the much simpler proof using the portrait map which was then also used in the final version of \cite{bib:StuVom11}. Corollary \ref{Cor: Vertices are the only lattice points} also implies that the set of standard imsets is exactly the set of all vertices of the standard imset polytope, again simplifying the lengthy proof from \cite{bib:StudVomHem10}. \eoproof

  Given a chordal undirected graph $G$, the corresponding characteristic imset $\chv_{G}$ is defined as the characteristic imset of any \DAG\, $\overrightarrow{G}$ Markov equivalent to $G$. The observation that characteristic imsets are unique representatives of Markov equivalence classes makes the definition correct.

  \begin{corollary}\label{Corollary: Entries of 1 for chordal graphs}
    Let $G$ be an undirected chordal graph over $N$. Then, for $T\subseteq N$, $|T|>1$, we have $\chv_{G}(T)=1$ if and only if $T$ is a clique in $G$.
  \end{corollary}

  \boproof As $G$ is the essential graph of an acyclic directed graph with no immoralities, we can direct the edges of $G$ in such a way that we obtain an equivalent acyclic directed graph $\overrightarrow{G}$ with no immoralities. To show the forward implication, let $T\subseteq N$ be given with $\chv_{G}(T)=1$. As $\chv_{\overrightarrow{G}}(T)=\chv_{G}(T)=1$, there exists some $i\in T$ such that $T\setminus\{ i\}\subseteq\pa_{\overrightarrow{G}}(i)$. Assume now, for a contradiction, that there are two nodes $j,k\in T\setminus\{i\}$ that are not connected by an edge in $\overrightarrow{G}$ (and hence $j$ and $k$ are not connected in $G$). Then, however, $j\rightarrow i\leftarrow k$ is an immorality in $\overrightarrow{G}$, a contradiction. Hence, all nodes in $T\setminus\{i\}$ must be pairwise connected by an edge in $G$. As they are all connected in $G$ by an edge to $i$, $T$ is a clique in $G$. To show the converse implication, let $T\subseteq N$ be a clique in $G$. Note that in $\overrightarrow{G}$, being an acyclic directed graph, the set $T$ must contain a node $i$ such that for all $j\in T$ the edge $\{i,j\}\in G$ is directed towards $i$ in $\overrightarrow{G}$. But then $T\setminus\{i\}\subseteq\pa_{\overrightarrow{G}}(i)$ and therefore, $\chv_{G}(T)=1$ by Theorem \ref{Theorem: Characteristic imsets are 0-1}. \eoproof

  Applying this observation to special undirected chordal graphs, namely to undirected forests, we obtain the following characterization.

  \begin{corollary}\label{Corollary: Characteristic imset for undirected forests}
    Let $G$ be an undirected forest having $N$ as the set of nodes. Then, for $T\subseteq N$, $|T|>1$, we have $\chv_{G}(T)=1$ if and only if $T$ is an edge of $G$, or, in other words,
    \[
      \chv_{G}=\left(\begin{array}{c}\chi(G) \\ \ve 0 \end{array}\right),
    \]
    where $\chi(G)$ denotes the characteristic vector of the edge-set of $G$.
  \end{corollary}

  Indeed, the only cliques of cardinality at least two in a forest are its edges. A similar result, in fact, holds for any acyclic directed graph $G$.

  \begin{corollary}\label{Corollary: Characteristic imset always contains chi(G)}
    Let $G$ be a \DAG\, over $N$ and $\bar G$ its underlying undirected graph. Then for any two-element subset $\{a,b\}\subseteq N$, we have $\chv_{G}(\{a,b\})=1$ if and only if $a\to b$ or $b\to a$ is an edge of $G$, or, in other words,
    \[
      \chv_{G}=\left(\begin{array}{c}\chi(\bar G) \\ * \end{array}\right),
    \]
    where $*$ denotes the remaining components of $\chv_{G}$.
  \end{corollary}

  \boproof This is an easy consequence of Theorem \ref{Theorem: Characteristic  imsets are 0-1}. If $\chv_{G}(T)=1$ for $T=\{ a,b\}$ then the only $i\in T$ with $T\setminus\{ i\}\subseteq\pa_{G}(i)$ are either $a$ or $b$. \eoproof

  Thus, $\chv_{G}$ is an extension of the characteristic vector $\chi(\bar G)$ of the edge-set of $\bar{G}$, which motivated our terminology. Let us now show how to convert $\chv_{G}$ back to the pattern graph $\pat(G)$ of $G$.

  \begin{theorem}\label{Theorem: Reconstruction of pattern graph from characteristic imset}
    Let $G$ be an acyclic directed graph over $N$ and $a,b\in N$ are distinct nodes. Then the following holds:
    \begin{itemize}
      \item[(1)] $a,b\in N$ are connected in $G$ iff $\chv_{G}(\{a,b\})=1$, otherwise $\chv_{G}(\{a,b\})=0$.
      \item[(2)] $a\to b$ belongs to an immorality in $G$ iff there exists some $i\in N\setminus\{a,b\}$ with $\chv_{G}(\{a,b,i\})=1$ and $\chv_{G}(\{a,i\})=0$. The latter condition implies $\chv_{G}(\{a,b\})=1$ and  $\chv_{G}(\{b,i\})=1$.
    \end{itemize}
  \end{theorem}

  \boproof The condition (1) follows from Corollary \ref{Corollary: Characteristic imset always contains chi(G)} and Theorem \ref{Theorem: Characteristic imsets are 0-1}.

   For (2) assume that $a\to b\leftarrow i$ is an immorality in $G$. Then $\chv_{G}(\{a,b,i\})=1$ by Theorem \ref{Theorem: Characteristic imsets are 0-1} and the necessity of the other conditions follows from (1). Conversely, provided that $\chv_{G}(\{a,b,i\})=1$, one of the three options $a\to i\leftarrow b$, $i\to a\leftarrow b$ and $a\to b\leftarrow i$ (with possible additional edges) occurs. Now, $\chv_{G}(\{a,i\})=0$ implies that $a$ and $i$ are not adjacent in $G$, which excludes the first two options and implies $a\to b\leftarrow i$ to be an immorality. \eoproof

  \begin{corollary}\label{Corollary: band determination}
    Let $G$ be a \DAG\ over $N$. The characteristic imset $\chv_{G}$ is determined uniquely by its values for sets of cardinality 2 and 3.
  \end{corollary}

  \boproof By Theorem \ref{Theorem: Reconstruction of pattern graph from characteristic imset} these values determine both the edges and immoralities in $G$. In particular, they determine the pattern $\pat(G)$. As explained in Section \ref{Ssec.Learning}, this uniquely determines the Bayesian network structure and, therefore, the respective standard and characteristic imsets. \eoproof

  More specifically, the components of $\chv_{G}$ for $|S|\geq 4$ can be derived iteratively from the components for $|S|\leq 3$ on the basis of the following lemma. A further simple consequence of the lemma below is that the entries for $|S|\geq 4$ are not linear functions of the entries for $|S|\leq 3$.

  \begin{lemma}\label{Lemma: band determination}
    Let $G$ be a \DAG\ over $N$, and $S\subseteq N$, $|S|\geq 4$. Then the following conditions are equivalent.
    \begin{itemize}
      \item[(a)] $\chv_{G}(S)=1$,
      \item[(b)] there exist $|S|-1$ subsets $T$ of $S$ with $|T|=|S|-1$ and $\chv_{G}(T)=1$,
      \item[(c)] there exist three subsets $T$ of $S$ with $|T|=|S|-1$ and $\chv_{G}(T)=1$.
    \end{itemize}
  \end{lemma}

  In the proof, by a {\em terminal node} within a set $T\subseteq N$ we mean $i\in T$ such that there is no $j\in T\setminus \{i\}$ with $i\to j$ in $G$.

  \boproof The implication $(a)\rightarrow (b)$ simply follows from Theorem \ref{Theorem: Characteristic imsets are 0-1}; $(b)\rightarrow (c)$ is trivial. To show $(c)\rightarrow (a)$ we first fix a terminal node $i$ within $S$. Now, $(c)$ implies there exist at least two sets $T\subseteq S$, $|T|=|S|-1$ which contain $i$. Let $\tilde{T}$ be one of them. Since $\chv_{G}(\tilde{T})=1$ by Theorem \ref{Theorem: Characteristic imsets are 0-1}, there exists $k\in\tilde{T}$ with $j\to k$ for every $j\in \tilde{T}\setminus\{k\}$. If $i\neq k$, then $i\to k$, which contradicts $i$ to be terminal in $S$. Thus, $i=k$. Since, those two sets $T$ cover $S$ one has $j\to i$ for every $j\in S\setminus \{i\}$ and Theorem \ref{Theorem: Characteristic imsets are 0-1} implies $\chv_{G}(S)=1$. \eoproof

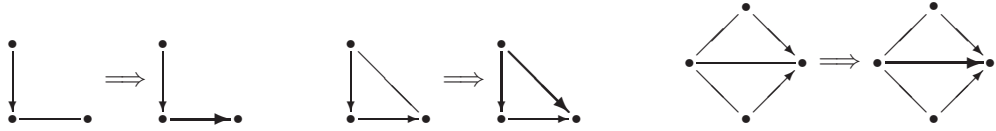
\begin{figure}
    \setlength{\unitlength}{1mm}
    \begin{center}
    \begin{picture}(128,20)
    \thinlines
    \put(2,2){\circle*{1}}
    \put(2,12){\circle*{1}}
    \put(12,2){\circle*{1}}
    \put(3,2){\line(1,0){8}}
    \put(2,11){\vector(0,-1){8}}
    \put(22,2){\circle*{1}}
    \put(22,12){\circle*{1}}
    \put(32,2){\circle*{1}}
    \put(22,11){\vector(0,-1){8}}
    \put(17,7){\makebox(0,0){$\Longrightarrow $}}
    \put(47,2){\circle*{1}}
    \put(47,12){\circle*{1}}
    \put(57,2){\circle*{1}}
    \put(48,2){\vector(1,0){8}}
    \put(47,11){\vector(0,-1){8}}
    \put(48,11){\line(1,-1){8}}
    \put(67,2){\circle*{1}}
    \put(67,12){\circle*{1}}
    \put(77,2){\circle*{1}}
    \put(68,2){\vector(1,0){8}}
    \put(67,11){\vector(0,-1){8}}
    \put(62,7){\makebox(0,0){$\Longrightarrow $}}
    \put(92,9.5){\circle*{1}}
    \put(99.5,2){\circle*{1}}
    \put(99.5,17){\circle*{1}}
    \put(107,9.5){\circle*{1}}
    \put(93,9.5){\line(1,0){13}}
    \put(93,10.5){\line(1,1){5.5}}
    \put(93,8.5){\line(1,-1){5.5}}
    \put(100.5,16){\vector(1,-1){5.5}}
    \put(100.5,3){\vector(1,1){5.5}}
    \put(117,9.5){\circle*{1}}
    \put(124.5,2){\circle*{1}}
    \put(124.5,17){\circle*{1}}
    \put(132,9.5){\circle*{1}}
    \put(118,10.5){\line(1,1){5.5}}
    \put(118,8.5){\line(1,-1){5.5}}
    \put(125.5,16){\vector(1,-1){5.5}}
    \put(125.5,3){\vector(1,1){5.5}}
    \put(112,9.5){\makebox(0,0){$\Longrightarrow $}}
    \thicklines
    \put(23,2){\vector(1,0){8}}
    \put(68,11){\vector(1,-1){8}}
    \put(118,9.5){\vector(1,0){13}}
    \end{picture}
    \end{center}
    \caption{Orientation rules for getting the essential graph.\label{fig.rules}}
  \end{figure}

  Theorem \ref{Theorem: Reconstruction of pattern graph from characteristic imset} allows us to reconstruct the essential graph for $G$. Indeed, the
  conditions (1) and (2) directly characterize the pattern graph $\pat(G)$. However, in general, there could be other arcs in the essential graph. Fortunately, there is a polynomial graphical algorithm transforming $\pat(G)$ into the corresponding essential graph $G^{*}$. More specifically, Theorem 3 in \cite{bib:Meek95} says that provided $\pat(G)$ is the pattern of an acyclic directed graph $G$ the repeated (exhaustive) application of the orientation rules from Figure \ref{fig.rules} gives the essential graph $G^{*}$.

  Finally, we wish to point out that Theorems \ref{Theorem: Characteristic imsets are 0-1} and \ref{Theorem: Reconstruction of pattern graph from characteristic imset} directly lead to a procedure for testing whether a given vector $\chv\in\Z^{2^{|N|}-|N|-1}$ is a characteristic imset for some (acyclic directed graph) $G$ over $N$. Using both theorems, one first constructs a candidate pattern graph, then a candidate essential graph, and then from it a candidate acyclic directed graph $G$. It remains to check whether the characteristic imset of $G$ coincides with the given vector $\chv$.

\section{Learning restricted Bayesian network structures}

  A lot of research is devoted to the topic of finding complexity results of the general problem of learning Bayesian network structures analyzing different optimization strategies, scoring functions and representations of data. For example, Chickering, Heckerman and Meek show the large-sample learning problem to be $\NP$-hard even when the distribution is perfectly Markovian \cite{bib:CHM04}. On the other side Chickering \cite{bib:Chi96} shows learning Bayesian network structures to be $\NP$-complete when using a certain Bayesian score. This remains valid even if the number of parents is limited to a constant.

\subsubsection*{Our assumptions}

  A reduction in complexity could be achieved by limiting the possible structures the Bayesian network can have. In the following, we will restrict our attention to learning \emph{decomposable models}, that is, learning the best \DAG\, among all \DAG s whose essential graphs are undirected (and thus also chordal). In fact, we assume that we are given an undirected graph $K$ over $N$ with an edge-set $\calE(K)$, not necessarily the complete graph, and we wish to learn a \DAG\ $G$ that maximizes the quality criterion and whose essential graph is an (undirected) subgraph of $K$ of a certain type. In particular, we are interested in learning undirected forests and spanning trees with and without degree bounds and in learning undirected chordal graphs.

  We wish to point out here that we make minimal assumptions on the database $D$ and on the quality criterion to be optimized. We only assume that the database $D$ over $N$ is \emph{complete}, that is, no data entry has a missing/unknown component (see Section \ref{Ssec.Learning}). For the quality criterion (= score function) we require that it is \emph{score equivalent} and \emph{decomposable}. In fact, instead of having $D$ and an explicit score function available, we only assume that we are given an evaluation oracle (depending on $D$) that, when queried on $G$, returns the value $\mathcal{Q}(G,D)$. Clearly, especially for larger databases $D$, computing a single score function value $\mathcal{Q}(G,D)$ may be expensive. By assuming a given evaluation oracle, we give constant costs to score function evaluations in our complexity results below.

  Finally, we wish to remind the reader that under our assumptions learning the best DAG representing $D$ becomes the problem of maximizing a certain {\em linear functional} (whose components depend on $D$) over the characteristic imsets (see Section \ref{Subsection: Algebraic approach to learning}). However, as this linear problem is in (exponential) dimension $2^{|N|}-|N|-1$, we cannot employ this transformation directly in our complexity treatment.

\subsection{Learning undirected forests and spanning trees}\label{Ssec:forests}

  By Corollary \ref{Corollary: Characteristic imset for undirected forests}, we know that every \DAG\ whose essential graph is an undirected forest $G$ has $\binom{\chi(G)}{\bf 0}$ as its characteristic imset. Thus, the problem of learning the best undirected forest is equivalent to maximizing a linear functional over such vectors $\binom{\chi(G)}{\bf 0}$ which in turn is equivalent to finding a maximum weight forest $G$ as a subgraph of $K$. The same argumentation holds for learning undirected spanning trees of $K$. These are two well-known combinatorial problems that can be solved in polynomial time via greedy-type algorithms (see e.g.\ \S\,40 in \cite{bib:Sch03}). We conclude the following statement.

  \begin{lemma}
    Given a node set $N$, an undirected graph $K=(N,\calE(K))$ and an evaluation oracle for computing $\mathcal{Q}(G,D)$. The problems of finding a maximum score subgraph of $K$ that is
    \begin{itemize}
      \item[(a)] a forest,
      \item[(b)] a spanning tree,
    \end{itemize}
    can be solved in time polynomial in $|N|$.
\end{lemma}

  Although $K$ is being part of the input, we need not state the complexity dependence with respect to the encoding length of $K$ explicitly here, since the encoding length $\langle K\rangle$ of $K$ is at least $|N|$. Moreover, we have $\langle K\rangle\in O(|N|^2)$.

  Chow and Liu \cite{bib:ChowLiu68} provided a polynomial time procedure (in $|N|$) for maximizing the maximum log-likelihood criterion which finds an optimal dependence tree (= a spanning tree). The core of their algorithm is the greedy algorithm and they apply it to a non-negative objective function. For their result, the complexity of computing the probabilities from data (and hence the objective/score function) is also omitted. A similar result was obtained by Heckerman, Geiger and Chickering \cite{bib:HGC95} for the Bayesian scoring criterion. Our result combines all of these previous results by only supposing a decomposable and score equivalent quality criterion.

  We wish to point out here that the well-known GES algorithm \cite{bib:Chi02,bib:Meek97}, which was designed to learn general Bayesian network structures, could be modified in a straight-forward way to learn undirected forests (among the subgraphs of $K$). Then the first phase of this new GES-type algorithm coincides with the greedy algorithm to find a maximum weight forest and the second phase of the algorithm cannot remove any edge. Thus, the modified GES algorithm always finds a best undirected forest (among the subgraphs of $K$) in time polynomial in $|N|$.

\subsection{Learning undirected forests and spanning trees with degree bounds}\label{Ssec:forests-bounded}

  Although the problems of learning undirected forests and of learning undirected spanning trees are solvable in polynomial time, learning an undirected forest/spanning tree with a given degree bound $\deg_G(i)\leq k<|N|-1$, $\forall i\in N$, is $\NP$-hard. For $k=1$ this problem is equal to the well-known problem of finding a maximum weight matching in $K$, which is in the general case polynomial time solvable (see \S\,30 in \cite{bib:Sch03}).

  \begin{theorem}\label{Theorem: LUFSTDB}
    Given a node set $N$, an undirected graph $K=(N,\calE(K))$ and an evaluation oracle for computing $\mathcal{Q}(G,D)$. Moreover, let $k\in\Z_+$ be a constant with $2\leq k<|N|-1$. Then the following statements hold.
    \begin{itemize}
      \item[(a)] The problem of finding a maximum score subgraph of $K$ that is a forest and that fulfils the degree bounds $\deg(i)\leq k$, $\forall i\in N$, is $\NP$-hard (in $|N|$) for any {\em fixed} (strictly) positive score function $\mathcal{Q}(.,D)$.
      \item[(b)] The problem of finding a maximum score spanning tree of $K$ that fulfils the degree bounds $\deg(i)\leq k$, $\forall i\in N$, is $\NP$-hard (in $|N|$) for any {\em fixed} score function $\mathcal{Q}(.,D)$.
    \end{itemize}
  \end{theorem}

  {\bf Remark.} Again, we have removed the explicit dependence on $\langle K\rangle$, since $\langle K\rangle\in O(|N|^2)$.

  \boproof We deduce part (b) from the following feasibility problem. In \S\,3.2.1 of \cite{bib:GarJohn79}, it has been shown that the task:

  \textsc{Bounded degree spanning tree}\\
  Instance: An undirected graph $K$ and a constant $2\leq k<|N|-1$\\
  Question: ``Is there a spanning tree for $K$ in which no node has degree exceeding $k$?''

  is $\NP$-complete by reduction onto the \textsc{Hamiltonian path problem}.

  Part (a) now follows by considering the subfamily of problems in which the linear objective takes only (strictly) positive values and, thus, every optimal forest (with the bounded degree) is a spanning tree. Hence, the problem of finding a maximum-weight forest (with a given degree bound) is equivalent to finding a maximum-weight spanning tree (with a given degree bound). As the feasibility problem for the latter is already $\NP$-complete, part (a) follows. \eoproof

  We wish to remark that Meek \cite{bib:Meek01} shows a similar hardness result for learning paths, i.e.\ spanning trees with upper degree bound $k=2$ for the maximum log-likelihood, the minimum description length and the Bayesian quality criteria.

\subsection{Learning chordal graphs}\label{Ssec:LCG}

  Undirected chordal graph models are the intersection of \DAG\ models and undirected graph models, known as Markov networks \cite{bib:Stud05}. In this section, we show that learning these models is $\NP$-hard.

  \begin{theorem}\label{Theorem: LCG}
    Given a node set $N$, an undirected graph $K=(N,\calE(K))$ and an evaluation oracle for computing $\mathcal{Q}(G,D)$. The problem of finding a maximum score chordal subgraph of $K$ is $\NP$-hard (in $|N|$).
  \end{theorem}

  \boproof We show that we can polynomially transform the following $\NP$-hard problem to learning undirected chordal graphs.

  \textsc{Clique of given size}\\
  Instance: An undirected graph $K$ and a constant $2\leq k\leq |N|-1$\\
  Question: ``Is there a clique set in $K$ of size at least $k$?''

  We define now a suitable learning problem that would solve this problem. By Corollary \ref{Corollary: Entries of 1 for chordal graphs}, we know that for any chordal graph $G$ the entry $\chv_{G}(T)$ is $1$ if and only if $T\subseteq N$, $|T|>1$ is a clique (otherwise this entry is $0$). Thus, the score function value for $G$ is determined by the values of the linear objective function ${\ve w}^\intercal\vex$ for the cliques $T$ in $G$. In particular, we can define the values for the cliques in such a way that when transforming the learning problem to the problem of maximizing ${\ve w}^\intercal\vex$ over the characteristic imset polytope (= the convex hull over all characteristic imsets) the entries ${\ve w}(T)$ are $0$ when $|T|<k$ and are positive when $|T|\geq k$. This implies that the maximum score among all chordal subgraphs of $K$ is positive iff there exists a chordal subgraph in $K$ containing a clique $T$ of size $|T|\geq k$. This happens iff $K$ has a clique of size at least $k$. \eoproof

\subsection{Learning chordal graphs with bounded size of cliques}

  Let us consider a variation of the previous task by introducing an upper bound $k$ for the size of cliques. If $k\leq 2$, we get the problems of learning undirected forests/matchings, which we already know are solvable in polynomial time (see Section \ref{Ssec:forests} and Section \ref{Ssec:forests-bounded}).

  For $k>2$, the corresponding problem is $\NP$-hard already for a fixed type of score function. This has been shown by Srebro \cite{bib:Sre01} for the maximized log-likelihood criterion (as a generalization of the work by Chow and Liu \cite{bib:ChowLiu68}).

\section{Conclusions}

  Let us summarize the main contributions of the paper. We introduced characteristic imsets as new simple representatives of Bayesian network structures, which are much closer to the graphical description. Actually, there is an easy transformation from the characteristic imset into the (essential) graph. Last but not least, the insight brought by the use of characteristic imsets makes it possible to offer elegant combinatorial proofs of (known and new) complexity results. The proofs avoid special assumptions on the form of the quality criterion besides the standard assumptions of score equivalence and decomposability.

  In our future work, we plan to apply these tools in the linear programming approach to learning. For this purpose we would like to find a general linear (facet-) description of the corresponding characteristic imset polytope or, at least, of a suitable polyhedral relaxation containing exactly the characteristic imsets as lattice points. Finding suitable polyhedral descriptions is also interesting and important for learning restricted families of Bayesian network structures, for example, for learning undirected chordal graphs.

  Finally, let us remark that a polyhedral approach to learning Bayesian network structures (using integer programming techniques) has been also suggested by Jaakkola {\em et.al.} \cite{bib:jaak10}, but their way of representing \DAG s is different from ours. Their representatives live in dimension $|N|\cdot 2^{|N|-1}$ and correspond to individual \DAG s, while ours live in dimension $2^{|N|}-|N|-1$ and correspond to Markov equivalence classes (of \DAG s).

\subsubsection*{Acknowledgements}

  The results of Milan Studen\'{y} have been supported by the grants GA\v{C}R n.\ 201/08/0539 and GAAV\v{C}R n.\ IAA100750603.

\end{document}